\documentclass[12pt,twoside,leqno]{article}
\usepackage{amsmath}
\usepackage{amssymb}
\usepackage{amsxtra}
\usepackage{amscd}
\usepackage{amsthm}
\usepackage[mathscr]{eucal}
\usepackage{hyperref}

\theoremstyle{plain}
\newtheorem{thm}[subsection]{Theorem}
\newtheorem{prop}[subsection]{Proposition}
\newtheorem{cor}[subsection]{Corollary}
\newtheorem{lem}[subsection]{Lemma}

\theoremstyle{definition}

\newtheorem{para}[subsection]{}

\newenvironment{pf}{\proof[\proofname]}{\endproof}
\begin{document}

\newcommand{\lr}[1]{\langle#1\rangle}
\newcommand{\ul}[1]{\underline{#1}}
\newcommand{\eq}[2]{\begin{equation}\label{#1}#2 \end{equation}}
\newcommand{\ml}[2]{\begin{multline}\label{#1}#2 \end{multline}}
\newcommand{\ga}[2]{\begin{gather}\label{#1}#2 \end{gather}}
\newcommand{\mc}{\mathcal}
\newcommand{\mb}{\mathbb}
\newcommand{\surj}{\twoheadrightarrow}
\newcommand{\inj}{\hookrightarrow}
\newcommand{\red}{{\rm red}}
\newcommand{\codim}{{\rm codim}}
\newcommand{\rank}{{\rm rank}}
\newcommand{\Pic}{{\rm Pic}}
\newcommand{\Div}{{\rm Div}}
\newcommand{\divi}{{\rm div}}
\newcommand{\Hom}{{\rm Hom}}
\newcommand{\Ext}{{\rm Ext}}
\newcommand{\im}{{\rm im}}
\newcommand{\fil}{{\rm fil}}
\newcommand{\gp}{{\rm gp}}
\newcommand{\Spec}{{\rm Spec}}
\newcommand{\Sing}{{\rm Sing}}
\newcommand{\Char}{{\rm char}}
\newcommand{\Tr}{{\rm Tr}}
\newcommand{\Gal}{{\rm Gal}}
\newcommand{\Min}{{\rm Min}}
\newcommand{\mult}{{\rm mult}}
\newcommand{\Max}{{\rm Max}}
\newcommand{\Alb}{{\rm Alb}}
\newcommand{\gr}{{\rm gr}}
\newcommand{\Ker}{{\rm Ker}}
\newcommand{\Lie}{{\rm Lie}}
\newcommand{\infi}{{\rm inf}}
\newcommand{\et}{{\rm \acute{e}t}}
\newcommand{\pole}{{\rm pole}}
\newcommand{\ti}{\times }
\newcommand{\modu}{{\rm mod}}
\newcommand{\Ab}[1]{{\mathcal A} {\mathit b}/#1}
\newcommand{\sA}{{\mathcal A}}
\newcommand{\sB}{{\mathcal B}}
\newcommand{\sC}{{\mathcal C}}
\newcommand{\sD}{{\mathcal D}}
\newcommand{\sE}{{\mathcal E}}
\newcommand{\sF}{{\mathcal F}}
\newcommand{\sG}{{\mathcal G}}
\newcommand{\sH}{{\mathcal H}}
\newcommand{\sI}{{\mathcal I}}
\newcommand{\sJ}{{\mathcal J}}
\newcommand{\sK}{{\mathcal K}}
\newcommand{\sL}{{\mathcal L}}
\newcommand{\sM}{{\mathcal M}}
\newcommand{\sN}{{\mathcal N}}
\newcommand{\sO}{{\mathcal O}}
\newcommand{\sP}{{\mathcal P}}
\newcommand{\sQ}{{\mathcal Q}}
\newcommand{\sR}{{\mathcal R}}
\newcommand{\sS}{{\mathcal S}}
\newcommand{\sT}{{\mathcal T}}
\newcommand{\sU}{{\mathcal U}}
\newcommand{\sV}{{\mathcal V}}
\newcommand{\sW}{{\mathcal W}}
\newcommand{\sX}{{\mathcal X}}
\newcommand{\sY}{{\mathcal Y}}
\newcommand{\sZ}{{\mathcal Z}}
\newcommand{\A}{{\mathbb A}}
\newcommand{\B}{{\mathbb B}}
\newcommand{\C}{{\mathbb C}}
\newcommand{\D}{{\mathbb D}}
\newcommand{\E}{{\mathbb E}}
\newcommand{\F}{{\mathbb F}}
\newcommand{\G}{{\mathbb G}}
\renewcommand{\H}{{\mathbb H}}
\newcommand{\I}{{\mathbb I}}
\newcommand{\J}{{\mathbb J}}
\newcommand{\M}{{\mathbb M}}
\newcommand{\N}{{\mathbb N}}
\renewcommand{\P}{{\mathbb P}}
\newcommand{\Q}{{\mathbb Q}}
\newcommand{\R}{{\mathbb R}}
\newcommand{\T}{{\mathbb T}}
\newcommand{\U}{{\mathbb U}}
\newcommand{\V}{{\mathbb V}}
\newcommand{\W}{{\mathbb W}}
\newcommand{\X}{{\mathbb X}}
\newcommand{\Y}{{\mathbb Y}}
\newcommand{\Z}{{\mathbb Z}}
\newcommand{\pic}{{\text{Pic}(C,\sD)[E,\nabla]}}
\newcommand{\ocd}{{\Omega^1_C\{\sD\}}}
\newcommand{\oc}{{\Omega^1_C}}
\newcommand{\al}{{\alpha}}
\newcommand{\be}{{\beta}}
\newcommand{\ta}{{\theta}}
\newcommand{\ve}{{\varepsilon}}
\newcommand{\phe}{{\varphi}}
\newcommand{\om}{{\overline M}}
\newcommand{\sym}{{\text{Sym}(\om)}}
\newcommand{\an}{{\text{an}}}
\newcommand{\bs}{{\backslash}}
\newcommand{\lra}{\longrightarrow}


\centerline{ } 
\vspace{45pt} 
\centerline{\LARGE Modulus of a rational map} 
\vspace{3pt}
\centerline{\LARGE into a commutative algebraic group}
\vspace{20pt}
\centerline{\large Kazuya Kato and Henrik Russell} 
\vspace{10pt}
\centerline{\large March 2010} 
\vspace{10pt}

\bigskip
\section{Introduction}
Let $X$ be a normal algebraic variety over a perfect field $k$, 
let $G$ be a commutative algebraic group over $k$, 
and let $\phe:  X\to G$ be a rational map. 
In this paper, we give a geometric definition of a modulus of $\phe$ 
as an effective divisor $\sum_v m(v) v$ on $X$. 
Here $v$ ranges over all codimension 1 points of $X$ at which 
$\phe$ is not defined as a morphism and $m(v)$ is a certain integer $\geq 1$. 
In the curve case, this definition coincides with Serre's definition \cite{Se}, 
which is based on the theory of local symbols. 
The case $k$ is of characteristic $0$ was explained in our previous paper 
\cite[\S5]{KR}. 
We discuss the positive characteristic case in this paper. 
We study properties of this modulus. 

An alternative way to define the modulus of $\phe$ is by using K-theoretic id\'ele class groups developed by the first author and Shuji Saito in \cite{KS}, 
as done in \cite{Ons} for surfaces. 
The coincidence of these two approaches follows from Prop.~\ref{hlf3}. 

This paper is related to the theory of
generalized Albanese varieties developed in the papers \cite{Ru}, \cite{Ru2} of the second author. In particular, the following fact will be proved in \cite{Ru2} by using this paper. If $X$ is proper smooth and if $Y$ is an effective divisor on $X$, $\phe$ factors through the generalized Albanese variety $\Alb(X, Y)$ of $X$ with modulus $Y$ if and only if (modulus of $\phe$) $\leq Y$. In the case $k$ is of characteristic $0$, this was proved in \cite[\S5]{KR} as a consequence of the theory in \cite{Ru}. 

The definition of modulus of $\phe$ is given in \S3 assuming Thm.~\ref{thmA}. The proof of this theorem is completed in \S5. In \S6 and \S7, we consider the relation of modulus with local symbols. In \S8, we consider the relation of modulus with field extensions. 

\section{Filtrations on additive groups and on Witt vector groups}

Let $K$ be a discrete valuation field, and let $O_K$ be the valuation ring of $K$.

\begin{para}\label{add}
For $m\geq 0$, we define 
$$\fil_m(K) =\{f\in K\;|\;v_K(f) \geq -m\}.$$
Here $v_K$ denotes the normalized  valuation of $K$. 

\end{para}

\begin{para}\label{pWitt} Let $p$ be a prime number and assume $K$ is of characteristic $p$. 
Let $W_n(K)$ be the set of Witt vectors of length $n$ with entries in $K$. For $m\geq 0$, define
$$\fil_mW_n(K)=\{(f_{n-1}, \dots, f_0)\;|\;f_j \in K, \; p^jv_K(f_j)\geq
-m \;\;(0\leq j \leq n-1)\}\subset W_n(K).$$
This filtration appeared in the paper \cite{Br} of Brylinski. In the case $n=1$, this filtration coincides with the filtration on $K=W_1(K)$  in \ref{add}.

 Let $F: W_n(K)\to W_n(K)$ be the map $(a_{n-1},\dots, a_0) \mapsto (a_{n-1}^p, \dots, a_0^p)$. For $m\in \N$, let  $$\fil^F_mW_n(K)=\sum_{j\geq 0} F^j(\fil_mW_n(K)) \subset W_n(K).$$

We have $\fil_0W_n(K)=\fil_0^FW_n(K)= W_n(O_K)$.

If we regard $W_n(K)$ as a subgroup of $W_{n+1}(K)$ via $V: W_n(K)\to W_{n+1}(K)$\; ; \;$(a_{n-1}, \dots, a_0)\mapsto (0, a_{n-1}, \dots, a_0)$, we have
$$\fil_mW_{n+1}(K)\cap W_n(K)=\fil_mW_n(K), \quad \fil^F_mW_{n+1}(K)\cap W_n(K)=\fil^F_mW_n(K).$$
\end{para}

\section{Modulus}

\begin{para} Let $X$ be a normal algebraic variety over a perfect field $k$. We regard $X$ as a scheme. 
Let $G$ be a commutative smooth connected algebraic group
over $k$, and let $\phe:X\to G$ be a rational map. We define the modulus $$\modu(\phe)=\sum_v \modu_v(\phe)v$$ of $\phe$ as an effective divisor on $X$, 
where $v$ ranges over all points of $X$ of codimension one and $\modu_v(\phe)\in \N$ is as follows. 

The case $k$ is of characteristic $0$ is already explained in \cite{KR}. (In \cite{KR}, we assumed that $X$ is proper smooth over $k$, but this condition is not used in the definition.) 

\end{para}

\begin{para}
First assume $k$ is algebraically closed. 

Let $0\to L\to G\to A \to 0$ be the canonical exact sequence of commutative algebraic groups, where $A$ is an abelian variety and $L$ is an affine smooth connected algebraic group. Write $L=L_m\times L_u$ where $L_m$ is multiplicative and $L_u$ is unipotent. Then since $k$ is algebraically closed, $L_m \cong (\G_m)^t$ for some $t\geq 0$. If $k$ is of characteristic $0$, $L_u\cong (\G_a)^s$ for some $s\geq 0$. Fix such isomorphism. If $k$ is of characteristic $p>0$, $L_u$ is embedded into a finite direct sum $\oplus_{i=1}^s W_{n_i}$ of Witt vector groups for some $s\geq 0$ and for some $n_i\geq 1$. Fix such embedding.

Let $K$ be the function field of $X$. Since it holds 
$H^1_{\text{fppf}}(\Spec(\sO_{X, x}), \G_m)=0$ and $H^1_{\text{fppf}}(\Spec(\sO_{X,x}), L_u)=0$ for any point $x$ of $X$, we have exact sequences
$$0 \to L(K) \to G(K) \to A(K) \to 0, \quad 0\to L(\sO_{X,x})\to G(\sO_{X,x}) \to A(\sO_{X,x})\to 0.$$
If $v$ is a point of $X$ of codimension one, since $A$ is proper and $\sO_{X,v}$ is a discrete valuation ring, we have $A(K)=A(\sO_{X,v})$. Hence the canonical map $L(K)/L(\sO_{X,v})\to G(K)/G(\sO_{X,v})$ is bijective. Take an element $l\in L(K)$ whose image in $G(K)/G(\sO_{X,v})$ coincides with the class of $\phe\in G(K)$.
In the case $k$ is of characteristic $0$, let $(l_i)_{1\leq i\leq s}$ be the image of $l$ in $(\G_a(K))^s=K^s$ under $L\to L_u\cong (\G_a)^s$. In the case $k$ is of characteristic $p>0$, let $(l_i)_{1\leq i\leq s}$  be the image of $l$ in $\oplus_{i=1}^s W_{n_i}(K)$ under  $L\to L_u \overset{\subset}\to \oplus_{i=1}^s W_{n_i}$. 

If $\phe \in G(\sO_{X,v})$, then we define $\modu_v(\phe)=0$. If $\phe \notin G(\sO_{X,v})$ and if the characteristic of $k$ is $0$ (resp. $p>0$), then we define $$\modu_v(\phe)=1+\max\{r(l_i)\;|\; 1\leq i\leq n\},\quad \text{where for $f\in K$ (resp. $W_{n_i}(K)$)},$$
$$r(f)=\min\{r\in \N\;|\; f\in \fil_r(K)\}\quad (\text{resp}.\;\; r(f)=\min\{r\in \N\;|\; f\in \fil^F_rW_{n_i}(K)\}).$$

In the case $k$ is of characteristic $0$, it is easy to see that $\modu_v(\phe)$ is independent of the choice of the isomorphism $L_u \cong (\G_a)^s$. In the case $k$ is of characteristic $p>0$, however, it is not so easy to prove
\end{para}

\begin{thm}\label{thmA} Let the notation be as above, and assume $k$ is of characteristic $p>0$. Then $\modu_v(\phe)$ is independent of the choice of the embedding $L_u \to \oplus_{i=1}
^s W_{n_i}$.

\end{thm} This theorem will be proved in \S5.

\begin{para} Now we do not assume $k$ is algebraically closed. Then by Galois descent for 
$\Gal(\bar k/k)$, we see that there is a unique effective divisor $\modu(\phe)$ on $X$ whose pull back to $X\otimes_k \bar k$ is the modulus of the rational map
$X\otimes_k \bar k \to G \otimes_k \bar k$.  
\end{para}

\section{Quotients of the filtrations}

Let $p$ be a prime number, and let $K$ be a discrete valuation field of characteristic $p$ with residue field $\kappa$. 

We study $\fil^F_mW_n(K)/\fil^F_{[m/p]}W_n(K)$ and its quotient $\fil^F_mW_n(K)/\fil^F_{m-1}W_n(K)$, for $m\geq 1$.  
Here for $x\in \R$, $[x]$ denotes $\max\{a\in \Z\;|\;a\leq x\}$ as usual.

\begin{prop}\label{directsum} {\rm(1)} The following sequence is exact.
$$0 \to  \oplus_{j\geq 0}\; \fil_{[m/p]}W_n(K) \overset{h}\to \oplus_{j\geq 0}\; \fil_mW_n(K) \to \fil^F_mW_n(K) \to 0,$$
where the third arrow is $(x_j)_j\mapsto \sum_j F^j(x_j)$, and $h$ is the map $(x_j)_j \mapsto (y_j)_j$ with $y_0 = F(x_0)$,  $y_j=F(x_j)-x_{j-1}$ for $j\geq 1$. 

\medskip

{\rm(2)} We have an isomorphism 
\begin{eqnarray*}
\oplus_{i\geq 0} \; \fil_mW_n(K)/(\fil_{[m/p]}W_n(K)+F(\fil_{[m/p]}W_n(K))) & \overset{\simeq} \to & \fil^F_mW_n(K)/\fil^F_{[m/p]}W_n(K) \; ; \\ 
(x_i)_i & \mapsto & \sum_i F^i(x). 
\end{eqnarray*}

\end{prop} 
\begin{pf} (1) We prove that for each $i\geq 0$, the sequence 
$$0 \to  \oplus_{j=0}^{i-1} \;\fil_{[m/p]}W_n(K) \overset{h_i}\to \oplus_{j= 0}^i \fil_mW_n(K) \to \sum_{j=0}^i F^j\fil_mW_n(K) \to 0$$
is exact, where $h_i$ is the restriction of $h$. We prove this by induction on $i$. The case $i=0$ is trivial. Assume $i\geq 1$. 
The non-trivial point is the exactness at the central term. Let $x=(x_j)_j$ be an element of $\oplus_{j= 0}^i \fil_mW_n(K)$ such that $\sum_j F^j(x_j)=0$. We prove that $x$ belongs to the image of $h_i$. 
We have $F^i(x_i)= -\sum_{j=0}^{i-1} F^j(x_j) \in \fil_{mp^{i-1}}W_n(K)$. Hence $x_i \in \fil_{[m/p]}W_n(K)$. Let $y=(y_j)_j$ be the element of $ \oplus_{j= 0}^{i-1}\; \fil_{[m/p]}W_n(K)$ defined by
$y_{i-1}=x_i$ and $y_j=0$ for $0\leq j< i-1$, and let $x'=x+h_i(y)$.  Then $x'\in \oplus_{j=0}^{i-1}\; \fil_mW_n(K)$. By induction on $i$, $(x'_j)_j$ is in the image of $h_i$. 

\medskip

(2) follows from (1) easily. \end{pf}

\begin{para}\label{delta} For a commutative ring $R$, let $\Omega_R^1=\Omega_{R/\Z}^1$ be the differential module of $R$. Then for any commutative ring $R$ over $\F_p$, there is a homomorphism 
$$\delta: W_n(R)\to \Omega^1_R\; ; \; (a_{n-1}, \dots, a_0)\mapsto \sum_i a_i^{p^i-1}da_i.$$ 

\end{para}

\begin{para} Let $\Omega^1_{O_K}(\log)$ be the differential module of $O_K$ with log poles defined by
$$\Omega^1_{O_K}(\log) = (\Omega^1_{O_K}\oplus (O_K \otimes_{\Z} K^\times))/N$$
where $N$ is the $O_K$-submodule of $\Omega^1_{O_K}\oplus (O_K \otimes_{\Z} K^\times)$ generated by $(da, -a \otimes a)$ for $a\in O_K-\{0\}$. We have canonical homomorphisms 
$\Omega^1_{O_K}\to \Omega^1_{O_K}(\log)$ and $K^\times \to \Omega_{O_K}(\log)\;;\;a \mapsto \text{class}(0, 1\otimes a)$. We denote the latter map by $d\log$.  If the condition

\medskip

(i) the completion of $K$ is separable over $K$

\medskip \noindent
is satisfied, then for a lifting $(b_i)_i$ of a $p$-base of $\kappa$ to $O_K$ and for a prime element $t$ of $K$, $\Omega^1_{O_K}$ (resp. 
 $\Omega_{O_K}^1(\log)$) is a free $O_K$-module with base $(db_i)_i$ and $dt$ (resp.  
$(db_i)_i$ and $d\log(t)$). 

(The condition (i) is equivalent to the condition that 
$(b_i)_i$ and $t$ form a $p$-base of $K$. Recall that for a field $F$ of characteristic $p$, a family $(b_i)_{i\in I}$ of elements of $F$ is called a $p$-base of $F$ if $F$ is generated over $F^p$ by $b_i$ ($i\in I$) as a field and for any subset $J$ of $I$ such that $J\neq I$, $F$ is not generated over $F^p$ by $b_j$ ($j\in J$). Recall also that if $(b_i)_i$ is a $p$-base of $F$, $(db_i)_i$ is a base of the $F$-module $\Omega_F^1$.)

Without the assumption (i), for any integer $j\geq 0$, $\Omega^1_{O_K} \otimes_{O_K} O_K/m_K^j$ (resp. 
$\Omega_{O_K}^1(\log)\otimes_{O_K} O_K/m_K^j$) is a free $O_K/m_K^j$-module with base $(db_i)_i$ and $dt$ (resp. $(db_i)_i$ and $d\log(t)$). This is because this group is invariant under the completion of $K$, and the condition (i) is satisfied of course if $K$ is complete.

\end{para}

\begin{prop}\label{deltam} For $m\geq 1$, the homomorphism $\delta$ (\ref{delta}) for $K$ induces an injective homomorphism
$$\delta_m : \fil_mW_n(K)/(\fil_{[m/p]}W_n(K)+F(\fil_{[m/p]}W_n(K))) \to \Omega^1_{O_K}(\log) \otimes_{O_K} m_K^{-m}/m_K^{-[m/p]}.$$

\end{prop}

\medskip
\begin{pf} The problem is the injectivity. By induction on $m$, it is reduced to the injectivity of
$$A:=\fil_mW_n(K)/(\fil_{m-1}W_n(K)+F(\fil_{[m/p]}W_n(K))) \to \Omega^1_{O_K}(\log) \otimes_{O_K} m_K^{-m}/m_K^{1-m}.$$ 
 We assume $K=\kappa((t))$ without a loss of generality.  Note that
$$\Omega_{O_K}^1(\log) \otimes_{O_K} m_K^{-m}/m_K^{1-m} \cong \Omega_\kappa^1 \oplus \kappa\;\;;\;\; $$ $$adb\otimes t^{-m}  \leftrightarrow (adb, 0) \;\;(a, b\in \kappa), \quad ad\log(t) \otimes t^{-m} \leftrightarrow (0, a)\;\;(a\in \kappa).$$

We define an increasing filtration $(A_i)_{-1\leq i\leq n-1}$ on $A$ as follows. For $-1\leq i\leq n-1$, let $A_i$ be the image of $\fil_mW_{i+1}(K)$ in $A$ under $V^{n-1-i} : W_{i+1}(K) \to W_n(K)$. Then as is easily seen, $A_i=A$ if $i\geq \text{ord}_p(m)$,  $A_{-1}=0$, and for $0\leq i\leq r:=\min(\text{ord}_p(m), n-1)$, we have an isomorphism 
$$\kappa \;\;(\text{resp}. \;\kappa/\kappa^p) \overset{\cong}\to A_i/A_{i-1}\quad \text{in the case $i=\text{ord}_p(m)$ (resp. $i<\text{ord}_p(m)$)};$$
$$a \mapsto (f_{n-1}, \dots, f_0)\;\text{with}\;f_j=at^{-mp^{-i}}\;\text{if}\;j=i, \;\;f_j=0\;\text{otherwise}.$$
If $a_i\in \kappa$ ($0\leq i\leq r$) and $f_i=a_it^{-mp^{-i}}$ for $0\leq i\leq r$ and $f_i=0$ for $r<i<n$, then
the image of $(f_{n-1}, \dots, f_0)\in \fil_mW_n(K)$ in $\Omega^1_{O_K}(\log) \otimes_{O_K} m_K^{-m}/m_K^{1-m}\cong \Omega^1_{\kappa} \oplus \kappa$
is 
$$ \Big(\sum_{i=0}^r a_i^{p^i-1}da_i,  -  \frac{m}{p^r}\cdot a_r^{p^r}\Big) \in \Omega^1_{\kappa} \oplus \kappa.$$  

For $i \geq 0$, 
let $B_i$ be the subgroup of $\Omega_\kappa^1$ 
generated by elements of the form $a^{p^j-1}da$
 with $a\in \kappa$ and $0 \leq j\leq i$. For example, 
 $B_0=d\kappa$. Let $B_{-1}=0$. The theory of Cartier isomorphism shows
 
 \medskip
 
 (1) 
  $\kappa/\kappa^p\overset{\simeq}\to B_i/B_{i-1}\;;\;a\mapsto a^{p^i-1}da$
  
  \medskip\noindent  for $i\geq 0$. For $0\leq i \leq r$, the image of the composition $A_i\to\Omega^1_\kappa\oplus \kappa\to \Omega^1_{\kappa}$ is contained in $B_i$, and the composition $\kappa/\kappa^p\overset{\simeq}\to A_i/A_{i-1}\to B_i/B_{i-1}$ is nothing but the isomorphism (1). If $\text{ord}_p(m)\leq n-1$ and $i=\text{ord}_p(m)$, the composition $A_i\to \Omega^1_\kappa\oplus \kappa\to \kappa$ kills $A_{i-1}$, and the composition $\kappa\overset{\cong}\to A_i/A_{i-1}\to \kappa$ coincides with the injective map $a\mapsto -m/p^r\cdot a^{p^r}$. This completes the proof of the injectivity in the proposition. \end{pf}

\begin{para} Let $O_K[F]$ be the non-commutative polynomial ring
defined by
$$O_K[F]=\Big\{\sum_{j\geq 0}  F^ja_j\;;\;a_j\in O_K\Big\}, \quad Fa=a^pF \;\;(a\in O_K).$$
For $m\in \N$, let
$$D_m = O_K[F] \otimes_{O_K} \Omega_{O_K}^1(\log) \otimes_{O_K} m_K^{-m}/m_K^{-[m/p]},$$
$$\bar D_m= \kappa[F]\otimes_\kappa (\Omega_{O_K}^1(\log) \otimes_{O_K} m_K^{-m}/m_K^{1-m}).$$

\end{para}

\begin{para} For $m\in \N$, by  Prop.~\ref{directsum} (2) and Prop.~\ref{deltam},
we have an injective homomorphism
$$\theta_m \;:\; \fil^F_mW_n(K)/\fil^F_{[m/p]}W_n(K) \to D_m(K)\;:$$
$$\sum_{j\geq 0}\; F^j(x_j) \mapsto \sum_j \; F^j\otimes \delta_m(x_j)$$
for $x\in \fil_mW_n(K)$. 

For $m\geq 1$, $\theta_m$ induces an injective homomorphism
$$\bar \theta_m: \fil^F_mW_n(K)/\fil^F_{m-1}W_n(K) \to \bar D_m.$$

\end{para}

\begin{para}\label{flat} For $m\geq 0$, we define a subgroup ${}^\flat\fil^F_mW_n(K)$ of $\fil^F_mW_n(K)$ 
as 
follows. 

Let ${}^\flat \bar D_m$ be the image of  $\kappa [F] \otimes_{\kappa} (\Omega^1_{O_K} \otimes_{O_K} m_K^{-m}/m_K^{1-m})$ (here we do not put a log pole) in $\bar D_m$. 
We have $${}^\flat \bar D_m\cong  \kappa [F] \otimes_{\kappa} \Omega^1_\kappa\otimes_{\kappa} m_K^{-m}/m_K^{1-m}.$$
Note that $$\bar D_m/{}^\flat \bar D_m\cong \kappa[F]\otimes_\kappa m_K^{-m}/m_K^{1-m};\quad F^ja  \otimes d\log(t) \otimes t^{-m} \leftrightarrow F^ja \otimes t^{-m}$$
where $a\in \kappa$ and $t$ is a prime element of $K$.

Let ${}^\flat\fil^F_mW_n(K)\subset \fil^F_mW_n(K)$ be the inverse image of ${}^\flat \bar D_m$ under 
$\bar \theta_m : \fil^F_mW_n(K) \to \bar D_m$. We have
$${}^\flat\fil^F_mW_n(K)= \sum_{j\geq 0} F^j({}^\flat\fil_mW_n(K))$$
where ${}^\flat\fil_mW_n(K)$ is the subgroup of $\fil_mW_n(K)$ consisting of all
elements $(f_{n-1}, \dots, f_0)$ which satisfy the following condition: 
 If the $p$-adic order $i$ of $m$ is $< n$, then $p^iv_K(f_i)> -m$. 

\medskip

We have injections
$$\fil^F_mW_n(K)/{}^\flat \fil^F_mW_n(K) \overset{\subset}\to \bar D_m/{}^\flat \bar D_m$$ 
$${}^\flat\fil^F_m W_n(K)/\fil^F_{m-1}W_n(K) \overset{\subset}\to {}^\flat \bar D_m$$
induced by $\bar \theta_m$. 

\medskip

As is easily seen, we have

\medskip

(1) For $m\geq 1$, ${}^\flat\fil^F_mW_n(K)\supset \fil^F_{m-1}W_n(K)$. If $m$ is prime to $p$, then ${}^\flat\fil^F_m W_n(K) = \fil^F_{m-1}W_n(K)$.

\medskip

(2) If $\kappa$ is perfect, then  ${}^\flat\fil^F_m W_n(K) = \fil^F_{m-1}W_n(K)$.

\end{para}

\medskip

\begin{para} The following relation with the refined Swan conductor in \cite{Ka2}, \cite{Ma} is proved easily. By Artin-Schreier-Witt theory, we have an isomorphism
$$W_n(K)/(F-1)W_n(K)\cong H^1(K, \Z/p^n\Z) := H^1(\Gal(K^{\text{sep}}/K), \Z/p^n\Z)$$
where $K^{\text{sep}}$ denotes the separable closure of $K$. As in \cite{Ka2}, let $\fil_mH^1(K, \Z/p^n\Z)$ be the image of $\fil_mW_n(K)$. 
\end{para}

\begin{prop} 
Let $\fil_mH^1(K, \Z/p^n\Z)\to \Omega_{O_K}^1(\log) \otimes_{O_K} m_K^{-m}/m_K^{1-m}$ ($m\geq 1$) be the refined Swan conductor in \cite{Ka2} whose kernel is $\fil_{m-1}H^1(K, \Z/p^n\Z)$. Then
 we have a commutative diagram
$$\begin{matrix}
\fil^F_mW_n(K) & \to & D_m/D_{m-1}= \kappa[F] \otimes_{\kappa} (\Omega_{O_K}^1(\log) \otimes_{O_K} m_K^{-m}/m_K^{1-m})\\
\downarrow & & \downarrow\\
\fil_mH^1(K, \Z/p^n\Z) & \to & \Omega_{O_K}^1(\log)\otimes_{O_K} m_K^{-m}/m_K^{1-m}.
\end{matrix}$$
Here the right vertical arrow is induced from the ring homomorphism $\kappa[F]\to \kappa\;;\;\sum_i F^ia_i \mapsto \sum_i a_i$ $\;(a_i\in \kappa)$. 
\end{prop}

\section{Homomorphisms and the filtrations}

Let $K$ be a discrete valuation field of characteristic $p>0$.

We assume here that we are given a perfect subfield $k$ of $O_K$.

\begin{para} Let $n, n'\geq 1$ and assume that we are given a homomorphism $h: W_n \to W_{n'}$ of algebraic groups  over $k$.
Let $h_1: \G_a \to \G_a$ be the homomorphism induced by $h$ on the subgroups $\G_a\subset W_n$ (embedded via $V^{n-1}$) and $\G_a\subset W_{n'}$ (embedded via $V^{n'-1}$). Since the endomorphism ring of $\G_a$ over $k$ is $k[F]$ where $F$ acts as $\G_a\to \G_a\;;\;x\mapsto x^p$, we can regard $h_1$ as an element of $k[F]$.

The following proposition is proved easily.

\end{para} 

\begin{prop}\label{prop1} 
{\rm(1)} The homomorphism $h$ sends $\fil^F_mW_n(K)$ into $\fil^F_mW_{n'}(K)$.

\medskip

{\rm(2)} We have a commutative diagram
$$\begin{matrix} \fil^F_mW_n(K) & \overset{\theta_m}\to& D_m(K)
\\
\downarrow && \downarrow\\
\fil^F_mW_{n'}(K) & \overset{\theta_m}\to & D_m(K)\end{matrix}$$
where the left vertical arrow is induced from $h$ and the right vertical arrow is the multiplication $x \mapsto h_1x$ by $h_1\in k[F]$.

\end{prop}
\begin{pf}
Homomorphisms $W_n\to W_{n'}$ are described by $F$, $V$, and the multiplication by elements of $W(k)$. For each of them, we can check easily that the proposition holds. \end{pf}

\begin{thm}\label{inj} Let $h: \oplus_{i=1}^s W_{n_i} \to \oplus_{j=1}^{s'} W_{n'_j}$ $(s, s'\geq 0, n_i, n'_j\geq 1)$ be an injective homomorphism defined over $k$. Let $m\geq 0$. Then for $x\in \oplus_{i=1}^s W_{n_i}(K)$, $x$ belongs to $\oplus_{i=1}^s \fil^F_mW_{n_i}(K)$ if and only if $h(x)$ belongs to $\oplus_{j =1}^{s'}\fil^F_mW_{n'_j}(K)$. 
\end{thm} 

\begin{pf} Let $h_1: \oplus_{i=1}^s \G_a\to \oplus_{j=1}^{s'} \G_a$ be the homomorphism induced by $h$ on the subgroups $\oplus_{i=1}^s \G_a\subset \oplus_{i=1}^s W_{n_i}$ and $\oplus_{j=1}^{s'} \G_a\subset \oplus_{j=1}^{s'} W_{n'_j}$. This $h_1$ is understood as a matrix with entries in $k[F]$. Since $h$ is injective, the homomorphism $$\Hom_\kappa(\G_a, \oplus_{i=1}^s \G_a)\to \Hom_\kappa(\G_a, \oplus_{j=1}^{s'} \G_a)\;;\;g \mapsto h_1\circ g$$
is injective, where $\Hom_\kappa$ means the set of homomorphisms of algebraic groups over $\kappa$. This means that the map $\oplus_{i=1}^s \kappa[F]\to \oplus_{j=1}^{s'} \kappa[F]\;;\;  x\mapsto h_1x$ is injective. Hence for $m\geq 1$, the map $\oplus_{i=1}^s \;\bar D_m \to \oplus_{j=1}^{s'}\; \bar D_m\;;\;x \mapsto h_1x$  is injective. By Prop.~\ref{prop1} (2), this proves that $h$ induces an injective homomorphism
$\oplus_{i=1}^s \;\fil^F_mW_{n_i}(K)/\fil^F_{m-1}W_{n_i}(K)\to \oplus_{j=1}^{s'} \;\fil^F_mW_{n'_j}(K)/\fil^F_{m-1}W_{n'_j}(K)$.
 \end{pf}

\medskip

\begin{para}\label{pf1} Proof of Thm.~\ref{thmA}. Let $Y=\oplus_i W_{n_i}$, Consider another embedding $L_u\to Y':=\oplus_{i'} W_{n_{i'}}$. Embed the push out $Y''$ of $Y \leftarrow L_u \to Y'$ into
a finite direct sum $Y''=\oplus_{i''} W_{n_{i''}}$. Then we have the third embedding $L_u \to  Y''$ and injective homomorphisms $Y\to Y''$ and $Y'\to Y''$ which are compatible with embeddings. By Thm.~\ref{inj}, $\modu_v(\phe)$ defined by the first (resp. second) embedding coincides with that defined by the third embedding.\qed
\end{para}

\section{Local symbols}

\begin{para}\label{curve} Let $k$ be an algebraically closed field, let $X$ be a normal algebraic curve over $k$, let $G$ be a commutative smooth connected algebraic group over $k$, and let $\phe: X\to G$ be a rational map. 
Then in \cite{Se}, the modulus of $\phe$ was defined  by using local symbols. We show that our definition of modulus coincides, in the curve case, with this classical definition. 

\end{para} 

\begin{para}\label{6.2} Let $k$, $X$, $G$ and $\phe$ be as in \ref{curve}, and let $K$ be the function field of $X$. For each point $v$ of $X$ of codimension one (that is, $v$ is a closed point of $X$), the local symbol map $$(\;,\;)_v : G(K) \times K^\times \to G(k)$$
is defined as in \cite{Se}. It is a $\Z$-bilinear map, and is continuous for the $v$-adic topology.  In \cite{Se}, the modulus of $\phe$ is defined as the right hand side of the equation in the following proposition.

\end{para}

\begin{prop} Let the notation be as in \ref{6.2}. Then 
 our $\modu_v(\phe)$ satisfies 
$$\modu_v(\phe)= \min\{m\in \N\;|\; (\phe, U^{(m)}_v)_v=0\}.$$
Here $U^{(m)}_v$ is the $m$-th unit group at $v$, that is, $U^{(m)}_v= \Ker(\sO_{X,v}^\times \to (\sO_{X,v}/m_{X, v}^m)^\times)$ where $m_{X, v}$ is the maximal ideal of $\sO_{X, v}$. 
\end{prop}

{\it Proof.} Let $0\to L \to G \to A \to 0$ be as in \S3. Since $(G(\sO_{X,v}), \sO_{X,v}^\times)_v$ vanishes and since $L(K)/L(\sO_{X,v})\to G(K)/G(\sO_{X,v})$ is bijective, we are reduced to the case $G=L$. If $k$ is of characteristic $0$, we are reduced to the cases $G=\G_m$ and $G=\G_a$. If $k$ is of characteristic $p>0$, by embedding $L_u$ to a finite direct sum of Witt vector groups as in \S3, we are reduced to the cases $G=\G_m$ and $G=W_n$. In the case $G=\G_m$, 
 the local symbol coincides with $(f, g) \mapsto (-1)^{v(f)v(g)}(g^{v(f)}/f^{v(g)})(v)$ where $v(?)$ denotes the $v$-adic normalized valuation and $(v)$ denotes the value at $v$. By using this fact, the case $G=\G_m$ is proved easily. In the case $G=\G_a$, the local symbol map is $(f, g) \mapsto \text{Res}(fd\log(g))$ where $\text{Res}$ is the residue map. By using this fact, in the case $k$ is of characteristic $0$, the case $G=\G_a$ is proved easily. In the case $k$ is of characteristic $p>0$ and $G=W_n$, it is sufficient to prove the following proposition.

\begin{prop}\label{perfect} Let $K=\kappa((t))$ with $\kappa$ a perfect field of characteristic $p>0$. For $m\geq 1$, let $U^{(m)}=1+t^m\kappa[[t]]\subset \kappa[[t]]^\times$. Let $(\;,\;)_K: W_n(K) \times K^\times \to W_n(\kappa)$ be the local symbol for $G=W_n$. 

\medskip

{\rm (1)} For $m\geq 0$, we have $(\fil^F_mW_n(K), U^{(m+1)}_K)_K=0$.  

\medskip
{\rm (2)} Let $\phe\in \fil^F_mW_n(K)$, and let $\sum_i F^ia_i$ be the image of $\phe$ under $\fil^F_mW_n(K)/\fil^F_{m-1}W_n(K)$ $\to$  $D_m/D_{m-1}\cong \kappa[F]$, where the last isomorphism is given by $F^ia\otimes d\log(t) \otimes t^{-m} \mapsto F^ia$ ($a\in \kappa$).
Then for $b\in \kappa$, the local symbol $(\phe, 1+bt^m)$ coincides with the image of $\sum_i (a_ib)^{p^{i+1-n}}\in \kappa $ under the injection $V^{n-1}: \kappa\to W_n(\kappa)$. 

\medskip

{\rm(3)} If $\kappa$ is an infinite field, then for any $m\geq 0$, we have $$\fil^F_mW_n(K)=\big\{\phe\in W_n(K)\;\big|\;
(\phe, U^{(m+1)}_K)_K=0\big\}.$$
\end{prop}

\begin{para} For the proof of Prop.~\ref{perfect}, we use the following explicit description of the local symbol map of $W_n$.

 Let $A=W_n(\kappa)[[t]][t^{-1}]$. We have the evident  surjective ring homomorphism $A\to K$, and an injective ring homomorphism
$$\phi_n: W_n(K) \to A\;;\; (a_{n-1}, \dots, a_0) \mapsto 
\sum_{0\leq i\leq n-1} p^{n-1-i}\tilde a_i^{p^i}.$$
Here $\tilde a_i$ is any lifting of $a_i$ to $A$. Note that
 $p^{n-1-i}\tilde a^{p^i}$ are independent of the choice of the  lifting. The differential module $\Omega_A^1$ is a free $A$-module of rank $1$ with basis $d\log(t)$. We have a well-defined homomorphism 
$$d\log: K^\times \to \Omega_A^1/pdA\;;\;a\mapsto d\log(\tilde a)$$
where $\tilde a$ denotes any lifting of $a$ to $A$. Let
$$\text{Res}: \Omega^1_A\to W_n(\kappa)\;;\;\sum_i a_it^id\log(t) \mapsto a_0.$$
Then the local symbol $(\;,\;)_K$ for $G=W_n$ is expressed as 
$$\text{(1)}\quad (f, g)_K=F^{1-n}\text{Res}(\phi_n(f)d\log(\tilde g))\quad \text{for $f\in W_n(K)$ and $g \in K^\times$}.$$
Here $F^{-1}: W_n(\kappa)\to W_n(\kappa)$ is the inverse map of $F: W_n(\kappa)\to W_n(\kappa)$. In the case $n=1$, this formula coincides with the formula
$(f, g)_K=\text{Res}(fd\log g)$ for $G=\G_a$.

By using the explicit formula (1) of the local symbol, we obtain (1) and (2) of Prop.~\ref{perfect}. (3) of Prop.~\ref{perfect} follows from (1) and (2) of Prop.~\ref{perfect}.

The authors are sure that the above formula (1) is written in some references, but they could not find. This (1) can be deduced from the formula (2) below.

Let $W_n\Omega_K^{\bullet}$ be the de Rham-Witt complex of $K$. Then $W_n\Omega^1_K$ is a $W_n(K)$-module and we have a homomorphism $d\log : K^\times \to W_n\Omega^1_K$. There is a residue map
$$\text{Res}: W_n\Omega_K^1\to W_n(\kappa)$$
(\cite[\S2]{Ka1}, \cite[\S2]{Rul}) which generalize the residue map $\Omega_K^1\to \kappa$ (the case $n=1$). By \cite[Chap.~III, Lem.~3]{KS}, we have
$$\text{(2)}\quad (f, g)_K=\text{Res}(f d\log(g))\quad \text{for}\;\; f\in W_n(K), \;g\in K^\times.$$
The above formula (1) follows from this formula (2) and from 
$$F^{1-n}\text{Res}(\phi_n(f)d\log(\tilde g))= \text{Res}(f d\log(g))\quad \text{for}\;\; f\in W_n(K), \;g\in K^\times.$$

\end{para}

\section{Higher dimensional local fields}

\begin{para} The above relation between modulus and local symbols for curves is generalized to the higher dimensional cases, by using local symbols for higher dimensional local fields define in \cite[Chap.~III]{KS}.

Let $p$ be a prime number, let $k_0$ be a perfect field of characteristic $p$, and define fields $k_r$ ($r\geq 1$) inductively by 
$$k_r=k_{r-1}((t_r)).$$ 

 Let $G$ be a commutative smooth connected algebraic group over $k_0$. Then the local symbol map
$$(\;,\;)_{k_r}\;: \;G(k_r) \times K^M_r(k_r) \to G(k_0)$$
is defined in \cite{KS} where $K^M_r$ denotes the $r$-th Milnor K-group.  

In the case $G=W_n$, this local symbol map is described as follows.
Define rings $A_r$ ($r\geq 0$) inductively by $A_0=W_n(k_0)$ and $A_r =A_{r-1}[[t_r]][t_r^{-1}]$ for $r\geq 1$. Then the local symbol map of $k_r$ for $W_n$ is described as
$$\text{(1)}\quad (f, g)_{k_r} =F^{1-n}\text{Res}(\phi_n(f)d\log(\tilde g))\quad \text{for $f\in W_n(k_r)$ and $g\in k_r^\times$} $$
where Res is the  map
$$\text{Res} : \Omega^r_{A_r}\to W_n(k_0)$$
defined to be the composition of the evident residue maps $\Omega^i_{A_i} \to \Omega^{i-1}_{A_{i-1}}$ ($1\leq i\leq r$) and $\phi_n: W_n(k_r)\to A_r$ is  defined in the same way as  $\phi_n$ in the previous paragraph, respectively. This (1) is deduced from the description of the local symbol map (\cite{KS})
$$(f, g)_{k_r} = \text{Res}(fd\log(g))\quad \text{for $f\in W_n(k_r)$ and $g\in k_r^\times$} $$
where Res is the residue map
$$\text{Res} : W_n\Omega_{k_r}^r\to W_n(k_0)$$
defined in \cite[\S2]{Ka1}. 
\end{para}

\begin{para}
By using the  explicit presentation (1) of the local symbol, we can obtain the following 
generalization Prop.~\ref{hdlf} of Prop.~\ref{perfect} to higher dimensional local fields. In Prop.~\ref{hdlf}, for $r\geq 1$, we show that the two filtrations $\fil^F_{\bullet} W_n(k_r)$ and ${}^\flat \fil^F_{\bullet} W_n(k_r)$ (which are defined with respect to the $t_r$-adic valuation of $k_r$) are related to certain two filtrations $U^{(\bullet)}_r$ and $V_r^{(\bullet)}$ on $K^M_r(k_r)$, respectively. 

 Fix $r\geq 1$. We define subgroups $U_r^{(m)}$ and $V_r^{(m)}$ of $K_r^M(k_r)$. 
For $m\geq 1$, let $U_r^{(m)}$ be the subgroup of $K_r^M(k_r)$ generated by all elements of the form $\{x, y_1, \dots, y_{r-1}\}$ such that $y_i\in k_r^\times$ and $x\in 1+t_r^mk_{r-1}[[t_r]]\subset k_{r-1}[[t_r]]^\times$. 
For $m\geq 0$, let $V_r^{(m)}$ be the subgroup of $K_r^M(k_r)$ generated by all elements of the form $\{x, y_1, \dots, y_{r-1}\}$ such that $y_i\in k_{r-1}[[t_r]]^\times$ and $x\in\Ker(k_{r-1}[[t_r]]^\times \to (k_{r-1}[t_r]/(t_r^m))^\times)$. Then
$$V_r^{(m-1)}\supset U_r^{(m)}\supset V_r^{(m)}\quad \text{for all $m\geq 1$}.$$
Let $U_r^{(0)}=V_r^{(0)}$. 

For $m\geq 1$, we have surjective homomorphisms
$$s_m: \Omega_{k_{r-1}}^{r-1}\to V_r^{(m)}/U_r^{(m+1)}\;;\;ad\log(b_1)\wedge \dots d\log(b_{r-1})\mapsto \{1+at_r^m, b_1, \dots, b_{r-1}\}$$
 $$s'_m: \Omega_{k_{r-1}}^{r-2}\to U_r^{(m)}/V_r^{(m)}\;;\;ad\log(b_1)\wedge \dots d\log(b_{r-2})\mapsto \{1+at_r^m, b_1, \dots, b_{r-2}, t_r\}$$
($a\in k_{r-1}$, $b_j\in k_{r-1}^\times$).
 
\end{para}

\begin{prop}\label{hdlf} Let $r\geq 1$. Define the filtrations $\fil^F_mW_n(k_r)$ and ${}^\flat\fil_mW_n(k_r)$ by using the $t_r$-adic discrete valuation of $k_r$.
Let $(\;,\;)_{k_r}: W_n(k_r)\times K^M_r(k_r)\to W_n(k_0)$ be the local symbol map 
 of $k_r$ for $G=W_n$. 

\medskip
{\rm (1)} For $m\geq 0$, we have $$\big(\fil^F_mW_n(k_r), U_r^{(m+1)}\big)_{k_r}=0, \quad
\big({}^\flat\fil^F_mW_n(k_r),  V_r^{(m)}\big)_{k_r}=0.$$

{\rm (2.1)} Let $m\geq 1$. $\phe\in \fil^F_mW_n(k_r)$, and let $\sum_i F^ia_i$ ($a_i\in k_{r-1}$) be the image of $\phe$ under $\fil^F_mW_n(k_r)/{}^\flat\fil^F_mW_n(k_r) \to \bar D_m/{}^\flat \bar D_m\cong k_{r-1}[F]$, where the last isomorphism is given by $F^ia \otimes d\log(t_r) \otimes t_r^{-m} \mapsto F^ia$ ($a\in k_{r-1}$). 
Then for $b\in \Omega_{k_{r-1}}^{r-1}$, the local symbol $(\phe, s_m(b))$ coincides with the image of $\sum_i (\text{Res}(a_ib))^{p^{i+1-n}}\in k_0$ under the injection $V^{n-1}: k_0\to W_n(k_0)$. Here Res is the residue map
$\Omega_{k_{r-1}}^{r-1}\to k_0$.  

\medskip
{\rm (2.2)} Let $m\geq 1$. Let $\phe\in {}^\flat\fil^F_mW_n(k_r)$, and let $\sum_i F^ia_i$ ($a_i\in \Omega^1_{k_{r-1}}$) be the image of $\phe$ under ${}^\flat\fil^F_mW_n(k_r)/\fil^F_{m-1}W_n(k_r) \to {}^\flat \bar D_m\cong k_{r-1}[F]\otimes_{k_{r-1}} \Omega^1_{k_{r-1}}$, where the last isomorphism is given by $F^ia \otimes w \otimes t_r^{-m} \mapsto F^ia\otimes w$ for $a\in k_{r-1}$, $w\in \Omega_{k_{r-1}}^1$.
Then for $b\in \Omega_{k_{r-1}}^{r-2}$, the local symbol $(\phe, s'_m(b))$ coincides with the image of $\sum_i (\text{Res}(a_i\wedge b))^{p^{i+1-n}}\in k_0 $ under the injection $V^{n-1}: k_0\to W_n(k_0)$. 

\medskip

{\rm(3)} If $k_0$ is an infinite field, then for any $m\geq 0$, we have $$\fil^F_mW_n(k_r)=\{\phe\in W_n(k_r)\;|\;
(\phe, U_r^{(m+1)})_{k_r}=0\},$$
$${}^\flat\fil^F_mW_n(k_r)=\{\phe\in W_n(k_r)\;|\;
(\phe, V_r^{(m)})_{k_r}=0\}.$$

\end{prop}

\begin{para}\label{hlf2} The following relation between modulus and higher dimensional local fields is deduced from  Prop.~\ref{hdlf}. Let $k$ be an algebraically closed field, let $X$ be a normal algebraic variety over $k$, let $G$ be a commutative smooth connected algebraic group over $k$, and let $\phe: X\to G$ be a rational map. 
Let $K$ be the function field of $X$. Let $v$ be a point of $X$ of codimension one. 

Let $r=\dim(X)$, let $k_0=k$, and define $k_i$ ($i\geq 1$) as above.  Assume $r\geq 1$ and assume that we are given a homomorphism of fields $K\overset{\subset}\to k_r$ such that $k_{r-1}[[t_r]]\cap K=\sO_{X, v}$, $t_r k_{r-1}[[t_r]]\cap K=m_{X, v}$, $k_{r-1}$ regarded as the residue field of $k_{r-1}[[t_r]]$ is separable over the residue field of $v$, and the ramification index of $k_{r-1}[[t_r]]$ over $\sO_{X,v}$ is $1$. (There are lot of such $K\to k_r$.)

\end{para}

\begin{prop}\label{hlf3} (1) For the local symbol map
$(\;,\;)_{k_r}: G(k_r) \times K^M_r(k_r)\to G(k)$, we have $$\modu_v(\phe)=\min\{m\in \N\;|\; (\phe, U_r^{(m)})_{k_r}=1\}$$
($1$ denotes the neutral element  of $G$). 

\medskip

(2) In the case $G=W_n$, if we endow $K$ with the discrete valuation associated to $v$, we have for any $m\geq 0$ $$\fil^F_mW_n(K)= \{f\in W_n(K)\;|\; (f, U_r^{(m+1)})_{k_r}=0\},$$
$${}^\flat\fil^F_mW_n(K) =\{f\in W_n(K)\;|\; (f, V_r^{(m)})_{k_r}=0\}.$$

\end{prop}

\section{Extension of local fields and the filtrations}

In this section, let $K$ be a discrete valuation field of characteristic $p>0$, and let $\kappa$ be the residue field of $K$.

We consider how  the filtrations $\fil^F_{\bullet} W_n(K)$ and ${}^\flat\fil^F_{\bullet} W_n(K)$ behave when the field  $K$  extends. 
In Thm.~\ref{thmB} and Thm.~\ref{thmC} below, we show how these filtrations are characterized by using extensions of $K$ with perfect residue fields.

The following lemma can be proved easily.
\begin{lem}\label{exten} 
Let $K'$ be a discrete valuation field containing $K$ such that $O_{K'}\cap K=O_K$ and $m_{K'}\cap K=m_K$. 
Let $m'=e(K'/K)m$ where $e(K'/K)$ is the ramification index of $K'$ over $K$.

\medskip

{\rm(1)} $\fil^F_mW_n(K) \subset \fil^F_{m'}W_n(K')$.

\medskip

{\rm(2)} For $m\geq 1$, we have a commutative diagram
$$\begin{matrix} \fil^F_mW_n(K) & \overset{\theta_m}\to& D_m(K)
\\
\downarrow && \downarrow\\
\fil^F_{m'}W_n(K) & \overset{\theta_{m'}}\to & D_{m'}(K').\end{matrix}$$

\end{lem}

\begin{cor}
Let $s_K(\phe) = \min\{m\in \N\;|\;\phe\in \fil^F_mW_n(K)\}$. Then $s_{K'}(\phe) \leq e(K'/K)s_K(\phe)$.
\end{cor}

\begin{cor}\label{extcor} Let $m\geq 1$.

\medskip

{\rm(1)}  The map $\fil^F_mW_n(K)/{}^\flat\fil^F_mW_n(K)\to \fil^F_{m'}W_n(K')/{}^\flat\fil^F_{m'}W_n(K')$ is injective if $e(K'/K)$ is prime to $p$, and is the zero map if $e(K'/K)$ is divisible by $p$.

\medskip

{\rm(2)} The map ${}^\flat\fil^F_mW_n(K)/\fil^F_{m-1}W_n(K) \to {}^\flat\fil^F_{m'}W_n(K')/\fil^F_{m'-1}W_n(K')$ is injective if  the residue field of $K'$ is separable over $\kappa$.

\end{cor}

\begin{cor}\label{sep}  
  In the case $e(K'/K)$ is prime to $p$ and the extension of the residue field in the extension $K'/K$ is separable, we have
$$s_{K'}(\phe)=e(K'/K)s_K(\phe).$$
\end{cor}

\begin{pf}
This follows from Cor.~\ref{extcor}.  \end{pf}

\begin{para}\label{perf} We consider what happens for extensions $K'$ of $K$, which have perfect residue fields. We consider the following $K'$.

\medskip

(1) $K'$ is a discrete valuation field containing $K$ such that $O_{K'}\cap K=O_K$ and $m_{K'}\cap K=m_K$, and such that the residue field of $K'$ is perfect.

\medskip

We also consider the following $K'$.

\medskip
(2) $K'$ is as in (1), but satisfies furthermore $e(K'/K)=1$.

\end{para}

\begin{thm}\label{thmB}
Let $\phe \in W_n(K)$. Then $$s_K(\phe)=\text{sup}\{e(K'/K)^{-1}s_{K'}(\phe)\;|\;K' \;\text{is as in \ref{perf} (1)}\}.$$
\end{thm}

For the filtration ${}^\flat\fil^F_mW_n(K)$, we have
\begin{thm}\label{thmC}
Let $\phe\in W_n(K)$. Then 
$$\min\{m \geq 1\;|\;\phe\in {}^\flat\fil^F_mW_n(K)\}= 1+\max\{s_{K'}(\phe)\;|\;\text{$K'$ is as in \ref{perf} (2)}\}.$$
\end{thm} 

We use the following lemma for the proofs of these theorems.

\begin{lem}\label{lemperf} Let $K'$ be as in \ref{perf} (1). 
Then for $m\geq 2$,  we have a commutative diagram with injective rows
 $$\begin{matrix} {}^\flat\fil^F_mW_n(K)/{}^\flat\fil^F_{m-1}W_n(K) & 
 \overset{\bar \theta_m}\longrightarrow & \kappa[F] \otimes_\kappa \Omega_\kappa^1\otimes_{\kappa} m_K^{-m}/m_K^{1-m})\\
 \downarrow & & \downarrow \\
 \fil^F_{em-1}W_n(K')/\fil^F_{em-2}W_n(K')  & 
 \overset{\bar \theta_{m-1}}\longrightarrow
 & (\kappa'[F]\otimes_{\kappa'} m_{K'}^{1-m}/m_{K'}^{2-m})/N $$ 
\end{matrix} $$
where $e=e(K'/K)$, $N= \kappa[F]\otimes_{\kappa} m_K^{1-m}/m_K^{2-m}$ if $e=1$, and $N=0$ if $e\geq 2$, and  the right vertical arrow is the map induced from
$$O_K[F]\otimes_{O_K} \Omega^1_{O_K}\otimes_{O_K} m_K^{-m}/m_K^{1-m} \to O_{K'}[F] \otimes_{O_{K'}}
\Omega^1_{O_{K'}}\otimes_{O_{K'}} m_{K'}^{-m}/m_{K'}^{1-m}.$$
\end{lem}

This is proved easily. 
\begin{para} 
The proofs of Thm.~\ref{thmB} and \ref{thmC}. 

For $K'$ as in \ref{perf} (1), ${}^\flat\fil^F_mW_n(K)\subset \fil^F_{e(K'/K)m-1}W_n(K')$ by \ref{flat} (2). Hence by Prop.~\ref{exten} (1),  
 it is sufficient to prove the following (1) and (2).

\medskip

(1) Let $m\geq 1$, and assume $\phe\in \fil^F_mW_n(K)$, $\phe\notin {}^\flat\fil^F_mW_n(K)$. Then for any $K'$ as in \ref{perf} (2), we have 
$s_K(\phe)=s_{K'}(\phe)$. 

\medskip

(2) Let $m\geq 2$, and assume $\phe\in \fil^F_{m-1}{}^\flat W_n(K)$, $\phe \notin \fil^F_{m-1}W_n(K)$. Then 
 for any integer $e\geq 1$, there is $K'$ as in 
\ref{perf} (1) such that $e=e(K'/K)$ and such that $s_{K'}(\phe)=em-1$.

\medskip

We prove (1) and (2).

(1) follows from Prop.~\ref{exten} (2) easily, by looking at the coefficient of $d\log(t) \otimes t^{-m}$ in the image of $\phe$ under $\bar\theta_m$. (Here $t$ denotes any prime element of $K$.)

We prove (2). 
Take a lifting  $(\tilde b_i)_{i\in I}$ of a $p$-base $(b_i)_{i\in I}$ of $\kappa$ to $O_K$. Let $$\kappa'=\cup_{r\geq 0}\; \kappa(T_i\;;\;i\in I)^{1/p^r},$$ 
where $T_i$ ($i\in I$) are indeterminates. Let $t$ be another indeterminate. Let 
$\pi$ be a prime element of $K$. 
 Then there is a unique homomorphism of fields 
 $K\to K':=\kappa'((t))$ which sends $O_K$ into $O_{K'}$, $m_K$ into $m_{K'}$, and sends $\tilde b_i$ ($i\in I$) to $b_i+T_it$ and $\pi$ to $t^e$. 
The right vertical arrow in the diagram in \ref{lemperf} sends $F^ja \otimes db_i$ ($a\in \kappa$) to $F^jaT_i$,  and sends $F^ja \otimes d\pi$ ($a\in \kappa$) to $F^ja$ if $e=1$ and to $0$ if $e\geq 2$. From this, we see that in the case $e=1$, the map
 $${}^\flat \fil^F_mW_n(K)/{}^\flat \fil^F_{m-1}W_n(K)\to \fil^F_{m-1}W_n(K')/\fil^F_{m-2}W_n(K')$$ is injective, and in the case $e\geq 2$, the map 
 $${}^\flat \fil^F_mW_n(K)/\fil^F_{m-1}W_n(K)\to \fil^F_{m-1}W_n(K')/\fil^F_{m-2}W_n(K')$$ 
is injective.  This proves (2). 
\end{para}

{\scshape
\begin{flushright}
\begin{tabular}{l}
Department of Mathematics \\  
University of Chicago \\ 
Chicago, IL 60637 \\ 
{\upshape e-mail: \texttt{kkato@math.uchicago.edu}}\\
{\upshape e-mail: \texttt{henrik.russell@uni-due.de}}\\
\end{tabular}
\end{flushright}
}

\end{document}